\documentclass{amsart}
\title{Locally isometric families of minimal surfaces}

\author{Aaron Peterson, Stephen Taylor}
\date{}
\usepackage{amsmath}
\usepackage{amssymb}
\usepackage{amsthm}
\usepackage{latexsym}
\usepackage{graphicx}
\usepackage{graphics}
\usepackage{url}

\numberwithin{equation}{section}

\theoremstyle{plain}
\newtheorem{thm}{\hspace{5pt}Theorem}
\theoremstyle{definition}
\newtheorem{coro}{\hspace{5pt} Corollary}
\theoremstyle{plain}

\def\vs{\vspace{7pt}}
\def\hs{\hspace{7pt}}

\def\hr{{}_1h}
\def\hi{{}_2h}
\def\gr{{}_1g}
\def\gi{{}_2g}

\begin{document}

\maketitle
\tableofcontents

\begin{abstract}

\noindent\textbf{Abstract.} We consider a surface $M$ immersed in $\mathbb{R}^3$ with induced metric $g=\psi\delta_2$ where $\delta_2$ is the two dimensional Euclidean metric.  We then
construct a system of partial differential equations that constrain $M$ to lift to a minimal surface via the Weierstrauss-Enneper representation demanding the
metric is of the above form. It is concluded that the associated surfaces connecting the prescribed minimal surface and its conjugate surface
satisfy the system.  Moreover, we find multiple symmetries of the PDE which each generate a one parameter family of surfaces isometric to a specified minimal
surface.

\end{abstract}

\section{Introduction}

\hs Given two $C^2$ functions $u$ and $v$ that satisfy Laplace's equation, a complex valued harmonic function $f$ is
defined via the combination: $f = u+iv$.  The Jacobian of such a function is given by $J_f=u_xv_y-u_yv_x$.
 We will only consider harmonic functions that are univalent (injective) with positive Jacobian on $\mathbb{D}=\{z:|z|<1\}$.  On a
 simply connected domain $D\subset\mathbb{C}$, a harmonic mapping $f$ has a canonical decomposition $f=h+\overline{g}$ where
 $h$ and $g$ are analytic in $D$, which is unique up to a constant.  The dilatation $\omega$ of a harmonic map $f$ is defined by
 $\omega\equiv g'/h'$.  The following theorem provides the link between harmonic univalent functions
and minimal surfaces:

\begin{thm}(Weierstrass-Enneper Representation).  Every regular minimal surface has locally an isothermal parametric
representation of the form

\begin{equation} \left(\mathrm{Re}\left\{\int^zp(1+q^2)dw\right\}, \mathrm{Im}\left\{\int^zp(1-q^2)dw\right\},
2\mathrm{Im}\left\{\int^zpqdw\right\}\right) \label{wrep}\end{equation}

in some domain $D\subset\mathbb{C}$, where $p$ is analytic and $q$ is meromorphic in $D$, with $p$ vanishing only at
the poles (if any) of $q$ and having a zero of precise order $2m$ wherever $q$ has a pole of order $m$.  Conversely,
each such pair of functions $p$ and $q$ analytic and meromorphic, respectively, in a simply connected domain $D$
generate through the formulas (\ref{wrep}) an isothermal parametric representation of a regular minimal surface.

\end{thm}

We will use (\ref{wrep}) in the following form:

\begin{coro} For a harmonic function $f=h+\overline{g}$, define the analytic functions $h$ and $g$ by $h=\int^zpd\zeta$
and $g=-\int^z pq^2d\zeta$.  Then the minimal surface representation
 (\ref{wrep}) becomes

\begin{equation}
\left(\mathrm{Re}\{h+g\},\mathrm{Im}\{h-g\},2\mathrm{Im}\left\{\int^z\sqrt{h'g'}d\zeta\right\}\right)\label{rep}
\end{equation}

\end{coro}

See \cite{Dur}, \cite{Clu} for a further introduction to harmonic mappings.  For an introduction to minimal surfaces see \cite{Die}.

\section{The isometric condition}

Let $\mathbf{x}(u,v)$ be a parametrization for a surface $M$ immersed in $\mathbb{R}^3$.  Set $z=u+iv$ and define $\phi=\partial\mathbf{x}/\partial z$.  Let $E$, $F$, and $G$ be
the coefficients of the metric induced in $\mathbb{R}^3$ by $\mathbf{x}(u,v)$.  We then have the relations

\begin{gather}
\phi^2 = \frac{1}{4}(E-G-2iF)\label{eq:met1}\\
\overline{\phi}^2=\frac{1}{4}(E-G+2iF)\label{eq:met2}\\
|\phi|^2=\frac{1}{4}(E+G)\label{eq:met3}
\end{gather}

where $\phi^2$ is notation for $\phi\cdot\phi$.  Inverting this system we find

$$E=\overline{\phi}^2+\phi^2+2|\phi|^2$$
$$F=i(\phi^2-\overline{\phi}^2)$$
$$G= -\overline{\phi}^2-\phi^2+2|\phi|^2$$

Since the Weierstrauss-Enneper representation theorem requires that $M$ has an isothermal parametrization, we require $E=G$ and $F=0$, which implies

$$\phi^2=0\qquad \overline{\phi}^2=0\qquad E=2|\phi|^2$$

The first two equations are identically satisfied.  Expanding the constraint on $E$, and using the identity

$$|\phi|^2=\frac{1}{4}|p|^2\left((1+q^2)(1+\overline{q}^2)+(1-q^2)(1-\overline{q}^2)+4q\overline{q}\right)$$

we find $E=|h'|^2+|g'|^2$.  Defining $\mathrm{Re}\{h\}={}_1h$,  $\mathrm{Im}\{h\}={}_2h$, $\mathrm{Re}\{g\}={}_1g$, and $\mathrm{Im}\{g\}={}_2 g$, we have the
Cauchy Riemann and isometric conditions:

\begin{equation} \hr_u-\hi_v=0\qquad \hr_v+\hi_u=0 \label{eqone}\end{equation}
\begin{equation} \gr_u-\gi_v=0\qquad \gr_v+\gi_u=0 \label{eqtwo}\end{equation}
\begin{equation} \hr_u^2+\hi_u^2+\gr_u^2+\gi_u^2-2\sqrt{(\hr_u^2+\hi_u^2)(\gr_u^2+\gi_u^2)}=0 \label{eqthree}\end{equation}

We now proceed to calculate the symmetry group of (\ref{eqone})-(\ref{eqthree}).  For an introduction to symmetry methods see \cite{Can}, \cite{Olv}, and \cite{Ste}.

\section{Symmetry analysis}

Preforming a symmetry analysis on a system of nonlinear equations is widely regarded as the best way to find exact solutions of
the system.  We will now summarize the Lie method in \cite{Olv}.

Consider a system of partial differential equations given by

$$\Delta_\nu(x,u^{(n)})=0,\quad \nu=1,\ldots,l$$

with $x=(x^1,\ldots,x^p)$ the set of independent variables and $u=(u^1,\ldots,u^q)$ the set of dependent variables where $1,\ldots, q$
run over the set of all partial derivatives of $u$ up to order $n$.  For $u=f(x)$, with $f:\mathbb{R}^p\rightarrow\mathbb{R}^q$ and
components $f^i$, $i=1\ldots q$, we define the $n$-th prolongation of $f$ to be

$$\mathrm{pr}^{(n)}f:\mathbb{R}^p\rightarrow U^{(n)}$$

given by $u^{(n)}=\mathrm{pr}^{(n)}f$, $u_J^a=\partial_Jf^a$ where $J$ is a multi-index running over the space of all possible
derivatives.  For example if we consider $u=f(x,y)$, we can compute

$$\mathrm{pr}^{(2)}f=(x,y)=(u;u_x,u_y;u_{xx},u_{xy},u_{yy})$$

The space $\mathbb{R}^p\times U^{(n)}$ is called the $n$-th order jet space of $\mathbb{R}\times U$.  The fundamental idea behind
the method of symmetry analysis is to view $\Delta_\nu$ as a map from the $n$-th order jet space into $\mathbb{R}^l$, and assuming
derivative terms occur as polynomials in the system, we can identify $\Delta_\nu$ with a subvariety in the jet space given by

$$\mathcal{L}_\Delta=\{(x,u^{(n)})|\Delta(x,u^{(n)})=0\}$$

Now let $M\subset\mathbb{R}^p\times U$ be open.  A symmetry group of $\Delta_{\nu}$ is a local group of transformations $G$ acting
on $M$ such that when $u=f(x)$ solves $\Delta_{\nu}$, then $u=g\cdot f(x)$ solves $\Delta_{\nu}$ for all $g\in G$ where defined.

Let $X$ be a vector field on $M$, and assume $X$ infinitesimally generates the symmetry ground $G$ of $\Delta_\nu$.  Then by projecting
$X$ into $M$ via the exponential map, we may construct a local one-parameter group $exp(\epsilon X)$.  We may then define the prolongation of $X$ as

$$\mathrm{pr}^{(n)}X=\frac{d}{d\epsilon}\mathrm{pr}^{(n)}[\exp(\epsilon X)](x,u^{(n)})\bigg|_{\epsilon=0}$$

where $\exp$ is the exponential map.  We also define the Jacobi matrix of $\Delta_\nu$ to be

$$J_{\Delta_\nu}(x,u^{(n)})=\left( \frac{\partial\Delta_\nu}{\partial x^i},\frac{\partial\Delta_\nu}{\partial u^a_J}  \right)$$

and say $\Delta_\nu$ is maximal if the rank of $J_{\Delta_{\nu}}$ is $l$.  The following theorem constrains the form of coefficients of
the $n$-th prolongation of an infinitesimal generator for the symmetry group.

\begin{thm}
Let

$$X=\xi^i(x,u)\frac{\partial}{\partial x^i}+\phi_a\frac{\partial}{\partial u^a}$$

then $X$ has prolongation

$$\mathrm{pr}^{(n)}X=X+\phi^J_a(x,u^{(n)})\frac{\partial}{\partial u^a_J}$$

where

$$\phi_a^J(x,u^{(n)})=D_J\left(\phi_a-\xi^au_i^a \right)+\xi^i \partial_i u^a_{J}$$

where subscripts on $u$ indicate partial derivatives and sums over repeated indices are implicit.

\end{thm}

The following may be called the fundamental theorem of the Lie method:

\begin{thm}
Let $\Delta_\nu$ be a system of differential equations of maximal rank.  If $G$ is a local group of transformations acting on $M$
and

$$\mathrm{pr}^{(n)}X[\Delta_\nu(x,u^{(n)})]=0$$

whenever $\Delta_\nu=0$, for every infinitesimal generator $X$ of $G$, then $G$ is a symmetry group of $\Delta_\nu$.

\end{thm}

We use these two theorems to calculate the coefficients of the infinitesimal generator.  We then exponentiate the infinitesimal
generator to obtain the symmetry group of the system.  Finally, we apply these symmetries to known and usually simple solutions
of the system to obtain new and hopefully more interesting solutions.

\vspace{7pt}

We now turn attention to (\ref{eqone})-(\ref{eqthree}). The infinitesimal generator of the above system is given by:

$$\mathbf{v}=c^u\partial_u+c^v\partial_v+c^{\hr}\partial_{\hr}+c^{\hi}\partial_{\hi}+c^{\gr}\partial_{\gr}+c^{\gi}\partial_{\gi}+c^{\psi}\partial_{\psi}.$$

Since the system is first order, we need only consider the first prolongation

$$\mathrm{pr}^{(1)}\mathbf{v}=\mathbf{v}+\hr^u\partial_{\hr_u}+ \hr^v\partial_{\hr_v}+\hi^u\partial_{\hi_u}+ \hi^v\partial_{\hi_v}$$
$$+\gr^u\partial_{\gr_u}+ \gr^v\partial_{\gr_v}+\gi^u\partial_{\gi_u}+ \gi^v\partial_{\gi_v}$$

where the $c^i$ are functions of $u,v,\psi,\hr,\hi,\gr$, and $\gi$.
Applying the first prolongation to the PDE system gives the following coefficient relationships:

$$\mathbf{v}_1=\partial_u\quad\mathbf{v}_2=\partial_v\quad\mathbf{v}_3=\partial_{\hr}\quad\mathbf{v}_4=\partial_{\hi}
\quad\mathbf{v}_5=\partial_{\gr}\quad\mathbf{v}_6=\partial_{\gi}$$
$$\mathbf{v}_7=-v\partial_u+u\partial_v\quad \mathbf{v}_8=-\hi\partial_{\hr}+\hr\partial_{\hi}\quad \mathbf{v}_9=-\gi\partial_{\gr}+\gr\partial_{\gi}$$
$$\mathbf{v}_{10} = u\partial_u+v\partial_v+\hr\partial_{\hr}+\hi\partial_{\hi}+\gr\partial_{\gr}+\gi\partial_{\gi}$$

Exponentiating these infinitesimal vector fields gives the symmetry transformations:

\begin{equation}  h\rightarrow h(z-s)\quad g\rightarrow g(z-s)\label{symthree}\end{equation}
\begin{equation}  h\rightarrow h(z-is)\quad g\rightarrow g(z-is) \label{symfour}\end{equation}
\begin{equation}  h\rightarrow h+s\quad g\rightarrow g\end{equation}
\begin{equation}  h\rightarrow h+is\quad g\rightarrow g\end{equation}
\begin{equation}  h\rightarrow h\quad g\rightarrow g+s\end{equation}
\begin{equation}  h\rightarrow h\quad g\rightarrow g+is\end{equation}
\begin{equation}  h\rightarrow h(e^{is}z)\quad g\rightarrow g(e^{is}z)\end{equation}
\begin{equation}  h\rightarrow e^{is}h\quad g\rightarrow g\end{equation}
\begin{equation}  h\rightarrow h\quad g\rightarrow e^{is}g\end{equation}
\begin{equation}  h\rightarrow e^s h(e^{-s} z)\quad g\rightarrow e^s g(e^{-s}z)\label{symlast}\end{equation}

In \cite{Bil} the minimal symmetry group for the real minimal surface
equation

$$u_{xx}(1+u_y^2)+u_{yy}(1+u_x^2)-2u_xu_yu_{xy}=0$$

was calculated.  Many of the translational symmetries and an $e^{s} f(e^{-s}x,e^{-s}y)$ symmetry were found.  We note that the analogue of $\mathbf{v}^{10}$
in \cite{Bil} is similar but different, since is constrains a Weierstrauss-Enneper representation of a surface and not a graph.

\section{Discussion and applications}

\begin{figure}[h]
\centerline{\hspace{0pt}\hbox{
\includegraphics[width=4.62236in,height=4.7in]{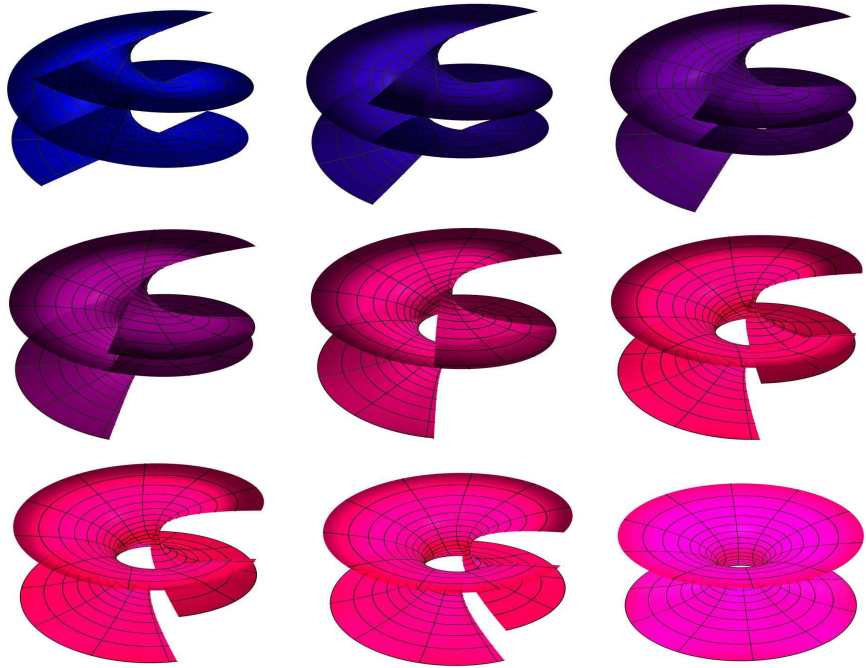}}}
\vspace{0pt}
\caption{\label{helcat} Helicoid to Catenoid Transformation.}
\end{figure}

\hspace{7pt} Consider the transformation $h\rightarrow e^{i\theta}h$, $g\rightarrow e^{i\theta}g$ which preserves the metric $E = |h'|^2 + |g'|^2$.  When
$\theta=0$, this is simply a minimal surface specified by defining $\psi$.  When $\theta=\pi/2$ we get the conjugate surface.  Thus all intermediate surfaces,
called associated surfaces, are isometric.
Since all minimal surfaces can be constructed from parts of a helicoid and catenoid \cite{Col}, the following examples are of interest.
First we draw attention to the catenoid,  given by $\psi=\cosh(v)^2$.  It's conjugate surface is the helicoid and associated surfaces between
the two are plotted over $\mathbb{D}$ in Figure \ref{helcat}.  Since all of the associated surfaces are isometric, geometrically they are equivalent.
However, note topologically the catenoid is $S^1\times\mathbb{R}$ where as the helicoid is $\mathbb{R}^2$.

\hspace{7pt} We now turn our attention to the other symmetries found in the analysis for the half catenoid.
We see that the symmetries generate surfaces that are topologically distinct from the catenoid, but geometrically identical as in the
above example.
Let $f$ be the harmonic mapping $f = h + \overline{g}$ where $$h = \frac{1}{2}\bigg(\frac{1}{2} \log \bigg[\frac{1 + z}{1 - z}\bigg] + \frac{z}{1 - z^2}\bigg)$$
and $$g = \frac{1}{2}\bigg(\frac{1}{2} \log \bigg[\frac{1 + z}{1 - z}\bigg] - \frac{z}{1 - z^2}\bigg)$$ which lifts to the catenoid.
We make the transformation in equation (\ref{symlast}) by letting $h \rightarrow e^sh(e^{-s}z)$ and $g \rightarrow e^sg(e^{-s}z)$. Figure 2 gives several plots
of this transformation for various $s$ values. The topology of the half catenoid is $\mathbb{R}^2$ for all $s$ up to
some value between $(0.3,0.4)$ where it changes to a punctured cylinder. Note as $s\rightarrow\infty$ that the minimal surfaces eventually degenerate to a line,
in a manner peculiarly similar to neckpinch singularities of the Ricci Flow.
\vs

\begin{figure}[h]
\centerline{\hspace{0pt}\hbox{
\includegraphics[width=6in,height=4in]{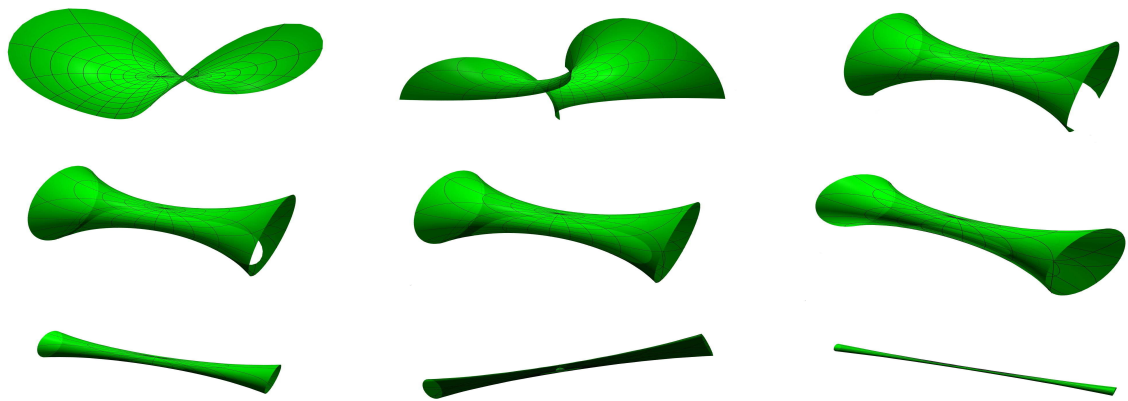}}}
\vspace{0pt}
\caption{\label{halfcat}Symmetry (\ref{symlast}) for $s = \{-1.2, -0.5, 0, 0.3, 0.4, 0.5, 1, 1.5, 3\}$.}
\end{figure}

When (\ref{symlast}) is applied to the helicoid, we find that the number of rotations of the helicoid about its axis are scaled.  Thus
we have:

\begin{thm} Let $S$ be the helicoid over $\mathbb{D}$ parameterized isothermally by\\
$\mathbf{x}=(\sinh u\sin v, \sinh u\cos v, -v)$. For helicoids $S_1$ given by $u\in(0,2\pi)$, $v\in(v_0,v_1)$, and $S_2$ by
$u\in(0,2\pi)$, $v\in(v_2,v_3)$ where $v_i\in\mathbb{R}$ then $S_1$ and $S_2$ are locally isometric. \end{thm}

\vs

It would be interesting to generalize the symmetry methods of this paper to higher dimensional Riemannian or Lorentizian manifolds.  One would
need a generalized Weierstrauss-Enneper representation for this task.  Moreover, we believe there are potential topological
theorems coming from symmetry (\ref{symlast}), which are connected to how the symmetry scales the domain of the graph under consideration.
For instance, if one calculates the one parameter family of minimal surfaces given by symmetry (\ref{symlast}) and a simply
connected minimal surface, does the topology always change from the plane to $S^1\times\mathbb{R}^1$ or some variant thereof?

\end{document}